\documentclass[a4paper,12pt]{article}

\usepackage{amsthm}
\usepackage{amsmath,amssymb,latexsym,amsfonts,mathrsfs}
\usepackage[dvipdfmx]{graphicx}

\usepackage{color}

\usepackage{fancyhdr}
\usepackage[top=30truemm, bottom=30truemm, left=25truemm, right=25truemm]{geometry}

\newcommand{\nn}{\nonumber}

\newcommand{\SCR}[1]{{\mathscr #1}}

\newcommand{\CAL}[1]{{\cal #1}
}

\newcommand{\MAT}[1]{\left(\begin{array}{cccccccccc}#1\end{array}\right)}

\theoremstyle{definition}
\newtheorem{Thm}{{\bf Theorem}}[section]

\newtheorem{Lem}[Thm]{{\bf Lemma}}

\newtheorem{Ass}[Thm]{{\bf Assumption}}

\newcounter{Exami}

\newcommand{\Proof}[2][Proof]{
\begin{proof}[{\bf #1}]
#2
\end{proof}
}

\begin{document}

\begin{flushleft}
{\bf \Large Strichartz estimates for Schr\"{o}dinger operators with square potential with time-dependent coefficients 
 } \\ \vspace{0.3cm} 
by {\bf \large Masaki Kawamoto } \\
Department of Engineering for Production, Graduate School of Science and Engineering, Ehime University, 3 Bunkyo-cho Matsuyama, Ehime, 790-8577,
Japan \\ 
Email: kawamoto.masaki.zs@ehime-u.ac.jp 
\end{flushleft}

\begin{center}
\begin{minipage}[c]{400pt}
{\bf Abstract}. {\small 
Strichartz estimates for a time-decaying harmonic oscillator were proven with some assumptions of coefficients for the time-decaying 
harmonic potentials. The main results of this paper are to remove these assumptions and to enable us to deal with the more general coefficient functions. Moreover, we also prove similar estimates for time-decaying homogeneous magnetic fields.
}
\end{minipage}
\end{center}

\begin{flushleft}
\end{flushleft}
\section{Introduction}

Let $n \in {\bf N}$. We consider the Hamiltonian 
\begin{align} \label{H0}
H_0 (t) = -\Delta /(2m) + \sigma (t) x^2/2, \qquad \mbox{on  } L^2({\bf R}^n),
\end{align}
where $x = (x_1,x_2,...,x_n) \in {\bf R} ^n $ and $m>0$ are the position and 
mass of a particle, respectively; $\Delta = \partial _1^2 + \cdots + \partial _n^2$ is the Laplacian.
We call $x^2= |x|^2= x_1^2+ \cdots + x_n^2$ the harmonic potential, and $\sigma (t)/2$ is a coefficient of the harmonic potential. 
In this paper, we impose the following conditions on the coefficient of the harmonic potential;
\begin{Ass} \label{A1}
Suppose $ \sigma \in L^{\infty} ({\bf R})$, and define $y_1 (t)$ and $y_2(t)$ as linearly independent solutions to 
\begin{align} \label{2}
y''_j(t) + \left( \frac{\sigma(t)}{m}\right) y_j(t) =0. 
\end{align}
Suppose that there exist constants $0 \leq \lambda <1/2$, $c_{1, \pm} \neq 0$ and $c_{2,\pm} \in {\bf R} $ such that 
\begin{align} \label{5}
\lim_{t \to \pm \infty} \frac{y_1(t)}{|t|^{1- \lambda}} = c_{1, \pm}, \quad 
\lim_{t \to \pm \infty} \frac{y_2(t)}{|t|^{\lambda}} = c_{2, \pm}.
\end{align} 
Moreover, $y_1 (t)$, $y_2 (t)$, $y_1' (t)$, and $y_2' (t)$ are continuous functions. 
\end{Ass}
The models of $\sigma (t)$, $y_1(t)$, and $y_2(t)$ with $\lambda \neq 0$ can be seen in, e.g., Geluk-Mari\'{c}-Tomi\'{c} \cite{GMT} 
(simplified models can be seen in Kawamoto \cite{K} and Kawamoto-Yoneyama \cite{KY}); the case of $\lambda = 0$ can be seen in, e.g., Naito \cite{Na} and Willett \cite{W}). 
Here, we remark that if only the conditions in this assumption are fulfilled, there is no need to assume $\sigma (t)$ is positive. 

To consider the Strichartz estimates, we define the Lebesgue space $L^{q,r}_{\lambda}$ as 
\begin{align*}
L_{\lambda}^{q,r} := \left\{ 
F \in \SCR{S}'({\bf R}^{n+1}) \, \Big| \, \left( 
\int_{- \infty}^{\infty} (1+t^2)^{- \lambda /2} \left\| 
F(t, \cdot)
\right\|^q_{r}  dt
\right)^{1/q}
 := \left\| 
F
\right\|_{q,r,\lambda} < \infty \right\} , 
\end{align*}
where we denote $\left\| \cdot \right\|_{L^{r}({\bf R}^n)} = \left\| \cdot 
\right\|_{r}$ throughout.

We define $U_0(t,s)$ as a propagator for $H_0(t)$, that is, a family of unitary operators $\left\{ U_0(t,s) \right\}_{(t,s) \in {\bf R}^2}$
so that each component satisfies 
\begin{align*}
& i \frac{\partial}{\partial t} U_0(t,s) = H_0(t) U_0(t,s), \quad i \frac{\partial}{\partial s} U_0(t,s) = -U_0(t,s) H_0(s), \\  
& U_0(t,s) U_0(s, \tau) = U_0(t, \tau).
\end{align*}
Moreover, we call a pair $(q,r)$ an admissible pair if it satisfies 
\begin{align}\label{adad1}
\frac{1}{q} + \frac{n}{2r} = \frac{n}{4}, \quad q >2, \quad r  \geq 2 .
\end{align}
Here we remark that the unique existence of the propagator $U_0(t,s)$ can be guaranteed by Theorem 1 in Fujiwara \cite{Fj}. The main theorem is the following: 
\begin{Thm} \label{T1}
Given the assumption \ref{A1}. Let $(q,r)$ and $(\tilde{q}, \tilde{r})$ be admissible pairs, $(q',r')$ and $(\tilde{q}', \tilde{r}')$ 
be a H\"{o}lder conjugate of $(q,r)$ and $(\tilde{q}, \tilde{r} )$, respectively. Then, for all  
$\phi \in L^{2}({\bf R}^n)$ and $F \in L_{\lambda}^{\tilde{q}',\tilde{r}'} $, the following Strichartz estimates hold: 
\begin{align*}
\left\| 
U_0(t,0) \phi 
\right\|_{q,r,\lambda} \leq C \left\| \phi  \right\|_{2}, 
\end{align*}
\begin{align*}
\left\| 
\int_0^t U_0(t,0) U_0(s,0)^{\ast} F(s) ds 
\right\|_{q,r, \lambda} \leq C \left\| 
F
\right\|_{\tilde{q}', \tilde{r}', - \lambda}
\end{align*}
and 
\begin{align*}
\left\| 
\int_{-t}^0 U_0(t,0) U_0(s,0)^{\ast} F(s) ds 
\right\|_{q,r, \lambda} \leq C \left\| 
F
\right\|_{\tilde{q}', \tilde{r}', - \lambda}.
\end{align*}
\end{Thm}
In the assumption of the admissible pairs, we removed the so-called end-points, that is, \eqref{adad1} with $q = 2$ and $n \geq 2 $. For the case of $\sigma (t) \neq 0$, we need to deal with the inhomogeneous dispersive estimate (see Lemma \ref{L1}). Then, the Hardy-Littlewood-Sobolev inequality yields the estimate 
\begin{align*}
\left\| 
\int_0^{|t|} U_0(t,0)U_0(s,0)^{\ast} F(s) ds 
\right\|_{q,r,0} \leq C \| F \|_{q',r',0} + C \| F \|_{q_{\lambda} ',r',0}, 
\end{align*} 
where $2/q_{\lambda} = n(1-\lambda) (1/2-1/r)$. With only this estimate, it is difficult to use the "duality argument", which plays a very important role in proving Strichartz estimates. Hence, in \cite{KY}, we avoid this difficulty by considering the time-weighted Lebesgue space $L_{\lambda}^{q,r}$, and we get 
\begin{align*}
\left\| 
\int_0^{|t|} U_0(t,0)U_0(s,0)^{\ast} F(s) ds 
\right\|_{q,r,\lambda} \leq C \| F \|_{q',r',-\lambda} , 
\end{align*}
which is quite suitable for using the duality argument. However, to prove 
\begin{align*}
 \| F \|_{q_{\lambda} ',r',0} \leq C  \| F \|_{q',r',-\lambda}
\end{align*}
we need to remove the end-point, even if $n \geq 3$. 

As usual for time independent selfadjoint energy $K$, $L^{2} ({\bf R}^n)$ can be divided by $L^{2} ({\bf R}^n)= L_p (K)  \oplus L_c (K)$ with the space of bound states; 
$$ L_p (K) := \left\{ \phi \in L^2({\bf R}^n) \, \middle| \, \lim_{r \to \infty} \sup_{t \in {\bf R}} \int_{|x| \geq r } \left| \left( e^{-itK} \phi \right)(x)\right|^2 dx = 0 \right\}
$$ 
and the space of scattering states;
$$ L_c (K) := \left\{ \phi \in L^2({\bf R}^n) \, \middle| \, \lim_{T \to \infty} \frac{
1}{T} \int_1^T \int_{|x| \leq r } \left| \left( e^{-itK} \phi \right) (x) \right|^2 dx dt = 0, \ {}^{\forall} r >0 \right\}
$$ 
Then it is well known that for all $\phi \in L_p(K)$, the global-in-time Strichartz estimates for $e^{-it K}$ fail since $ e^{-itK} \phi  $ is the non-dispersive solution, and in such case sometimes the local-in-time Strichartz estimate has been considered; e.g., Yajima-Zhang \cite{YZ}, and Mizutani \cite{Mi}, among others. Here we note that the space $L_p(K)$ can be characterized by eigenfunctions of $K$ with using the fact that $e^{-itK} \phi_{\lambda} = e^{-it \lambda} \phi_{\lambda}$, where $\phi_{\lambda}$ is an eigenfunction and $\lambda$ is an eigenvalue of $K$, respectively. On the other hand, if energy depends on time, to characterize bound state through eigenfunctions is difficult since the propagator can not be written like exponential function of energy. Indeed, for $H_0(t)$, under the assumption \ref{A1}, one can prove that for all $r>0$ and $\phi $ in some dense space on $L^2({\bf R}^n)$, the inequality 
\begin{align*}
\lim_{T \to \infty} \frac{1}{T} \int_{1}^T \int_{|x| \leq r} \left| U_0(t,0) \phi \right|^2 {dt} = 0
\end{align*}
holds, see e.g., Proposition 2.1 of Ishida-Kawamoto \cite{IK}, i.e., $H_0(t)$ has no bound states in spite of $H_0(t)$ has eigenvalues for any fixed $t$. 

As for $H_0(t)$, for all bounded $\sigma (t)$ (which includes $\sigma (t) = \sigma_0 >0$), the local-in-time dispersive estimate 
\begin{align*}
\left\| U_0(t,s) \phi  \right\| _{L^{\infty} ({\bf R}^n)} \leq C |t-s|^{-n/2} \|  \phi \|_{L^1({\bf R}^n)}, \quad |t-s| \leq {}^{\exists} R (>0)
\end{align*}
holds by Theorem 2 in \cite{Fj}, and as the direct consequence of \cite{KT}, the local-in-time Strichartz estimates also hold. Hence we focus on what kind of $\sigma (t)$ enable the global-in-time Strichartz estimates true. The most simplest model of $\sigma (t)$ is $\sigma (t) \equiv 0$ and in this case global-in-time Strichartz estimates were first considered by Strichartz \cite{S} for the case $\sigma (t) \equiv 0$, and the condition of admissible pairs stated in \cite{S} were relaxed by many researchers (e.g., Yajima \cite{Ya1}, Ginibre-Velo \cite{GV}, and Keel-Tao \cite{KT}). As the other model, we showed $\sigma (t)$ with \eqref{5} gives global-in-time Strichartz estimates; We expect this paper will be a first step to be found some other mathematical models of $\sigma (t)$ on which global-in-time Strichartz estimates hold. 

For time-dependent systems, recently dispersive and Strichartz estimates for time-dependent systems ($H(t) = -\Delta /(2m) + V(t,x)$ with $|V(t,x)| \to 0$ for $|x| \to \infty$)
have also been considered, e.g., Goldberg \cite{G}, Pierfelice \cite{P}, and Rodnianski-Schlag \cite{RS}. However, if the coefficient of harmonic potential depends on the time, the results of Strichartz estimates are very few. Carles \cite{Ca} studied a time-local Strichartz estimate for
the case of $\sigma (t) \neq 0$ with $\lambda =0$. In that paper \cite{Ca}, by using the lens
transform, it was revealed that the solution of $i \partial _t u(t,x) = H_0(t) u(t,x)$ can be written in the form 
\begin{align} \label{adad2}
u(t,x) =  e^{im \zeta _1'(t) x^2/(2 \zeta _1 (t) )} \left( |\zeta _1(t)|^{-n/2}  v \left( 
\frac{\zeta _2 (t)}{  \zeta_1(t)  } , \frac{x}{| \zeta _1(t) | }
\right) \right), 
\end{align} 
where $v(t,x) $ is a solution of $i \partial _t v(t,x) = (- \Delta /2m) v(t,x)$ with $v(0,x) =v_0$, and $\zeta _j (t)$, $j=1,2$ are the solutions to 
\begin{align*}
\zeta _j ''(t) + \frac{\sigma (t)}{m} \zeta _j(t) =0 , \quad 
\begin{cases}
\zeta _1(0) =1, \\ 
\zeta _1'(0) =0,
\end{cases}
\ 
\begin{cases}
\zeta _2(0) = 0, \\ 
\zeta _2'(0) =1, 
\end{cases}
\end{align*}
respectively. Here, by introducing the so-called generator of the dilation group $A = x \cdot (-i \nabla) + (-i \nabla ) \cdot x$, we get 
\begin{align*}
e^{-i \log |\zeta _1 (t)| A /2} f(t,x) = |\zeta _1 (t)|^{-n/2} f\left(t,  \frac{x}{|\zeta _1 (t)|} \right);
\end{align*}
see \S{2} in \cite{KY}. Hence, the equation \eqref{adad2} is equivalent to 
\begin{align*}
u(t,x) = e^{i m \zeta _1 '(t) x^2/(2 \zeta _1 (t))} e^{-i (\log |\zeta _1 (t) | ) A/2}e^{i \zeta _2 (t) \Delta /(2m \zeta _1(t))}  v_0. 
\end{align*} 
Using the commutator calculation, we get 
\begin{align*}
e^{-i (\log |\zeta _1 (t) | ) A/2}e^{i \zeta _2 (t) \Delta /(2m \zeta _1(t))} = e^{i \zeta _1(t)\zeta _2(t) \Delta /(2m)} e^{-i (\log | \zeta _1(t)| )A/2} , 
\end{align*} 
which provides 
\begin{align*}
u (t,x) = e^{i m \zeta _1 '(t) x^2/(2 \zeta _1 (t))} e^{i \zeta _1(t)\zeta _2(t) \Delta /(2m)} e^{-i (\log | \zeta _1(t)| )A/2} v_0.
\end{align*} 
This formula locally coincides with the following Korotyaev's factorization formula:
\begin{align} \label{1} 
U_0(t,0) = e^{ mi \zeta _1'(t) x^2/(2 \zeta _1(t))} e^{i \zeta _1(t)\zeta _2(t) \Delta /(2m)} e^{-i (\log | \zeta _1(t)| )(x \cdot (-i \nabla) + 
(-i \nabla ) \cdot x) /2} {\bf S}^{\nu (t)} ,  
\end{align} 
where $({\bf S} f) (x) = e^{-in\pi/2} f(-x)$ for $f\in L^2({\bf R}^n)$ and $\nu(t)$ is the number of zeros of elements of 
$\left\{ \tau \in [0,t] \, | \, \zeta _1(\tau) =0 \right\}$ for $t \geq 0$ or $\left\{ \tau \in [t,0] \, | \, \zeta _1(\tau) =0 \right\}$ for $t \leq 0$,
(see Korotyaev \cite{Ko} and \cite{KY}). Conversely, in \cite{KY}, using this factorization formula succeeded in deducing the dispersive estimates for $U_0(t,0)$ 
((1.6) in \cite{KY}), which enabled us to prove the Strichartz estimate.
However, in this case, the dispersive estimate, which is deduced from this factorization, can be written as 
\begin{align} \label{adad3}
\left\| 
U_0(t,0) U_0(s,0)^{\ast} \phi 
\right\|_{L^q({\bf R}^n)} & \leq C |\zeta _1 (t) \zeta _2 (s) - \zeta _1(s ) \zeta _2(t)|^{-n(1/2-1/q)} \left\| 
\phi
\right\|_{L^{q'} ({\bf R}^n)} , \\ 
 \Big( \mbox{not    }  &  \leq  C |\zeta _1(t-s)|^{-n(1/2-1/q)} \left\| 
\phi
\right\|_{L^{q'} ({\bf R}^n)} \Big) .\nn
\end{align} 
Here, we divide the region $(t,s) \in{\bf R}^2$ into $|t-s| \geq R$ and $|t-s| \leq R$ for some $R>0$. Then, it seems that we can apply the result of \cite{Fj} on the region $|t-s| \leq R$ and apply the result of \cite{Ya1}, \cite{GV}, \cite{KT} on the region $|t-s| \geq R$. However, dealing with these arguments is difficult because the coefficient on the right-hand-side of \eqref{adad3} has singularities, even if $|t-s| \geq R$ since the factorization \eqref{1} has singularities in the case of $\zeta _1(t) = 0$. Hence, in \cite{KY}, some specific assumptions on $\sigma (t)$ to 
deal with such difficulties were assumed. The aim of this paper is to remove these assumptions. The approach is to use the 
factorization formula given by \cite{AK}, which is different from that of \eqref{1}. 
Because this formula has no singularities for every $t \in {\bf R}$, we can deal with more general coefficient functions. 

Next, let us consider the Hamiltonian 
\begin{align}\label{Mg01}
H_{2,L} (t) = \frac{1}{2m} \left(-i \partial_1 + \frac{qB(t)}{2} x_2 \right)^2 +   \frac{1}{2m} \left(-i \partial_2 - \frac{qB(t)}{2} x_1 \right)^2, \quad 
\mbox{on  } L^2({\bf R}^2)
\end{align}
and 
\begin{align} \label{Mg02}
H_{3,L} (t) = \frac{1}{2m} \left(-i \partial_1 + \frac{qB(t)}{2} x_2 \right)^2 +   \frac{1}{2m} \left(-i \partial_2 - \frac{qB(t)}{2} x_1 \right)^2 - 
\frac{\partial _3^2}{2m}, \quad \mbox{on  } L^2({\bf R}^3),  
\end{align}
where $q \neq 0$ is the charge of a charged particle and $B(t) \in L^{\infty} ({\bf R})$ is the intensity of the time-dependent magnetic field 
${\bf B} (t) = (0,0, B(t))$. We define $U_{2,L} (t,0)$ and $U_{3,L} (t,0)$ as propagators for $H_{2,L} (t)$ and $H_{3,L} (t)$, respectively. Here, we note 
\begin{align*}
H_{2,L} (t) = -\frac{1}{2m} \Delta _{(2)} +  \frac{q^2B(t)^2}{8m}  |x_{(2)}|^2 - \frac{qB(t)}{2m}L_{(2)}, 
\end{align*}
where $\Delta _{(2)} = \partial _1^2 + \partial _2^2$, $|x_{(2)}|^2 = x_1^2+x_2^2$ and $L_{(2)}= x_1 (- i\partial _2)- x_2 (-i \partial _1)$; $L_{(2)}$ is called the angular momentum. The unique existence of propagators $U_{j,L} (t,0)$, $j=2,3$ can be guaranteed by Theorem 1 in Yajima \cite{Ya2}. 

If $B(t)$ is a constant $B \neq 0$, $H_{2,L}$ is called the Landau Hamiltonian and it is well known that 
the solution of $i \partial _t \phi (t,x) = H_{2,L} \phi (t,x)$ has the refocusing property; it also well known that in this case 
time-global Strichartz estimates fail.
Hence, the previous studies of Strichartz estimates for the magnetic Schr\"{o}dinger operator consider the case where the magnetic potential is decaying in $x$ or 
is localized near the origin; see, e.g., D'Ancona-Fanelli-Vega-Visciglia \cite{DFVV}, Erdo\v{g}an-Goldberg-Schlag \cite{EGS}, and  Stefanov \cite{S}. Such magnetic fields 
are inhomogeneous. For the homogeneous case, it seems necessary to assume that the intensity of the magnetic field depends on time.

Here, we rewrite $q^2B(t)^2/(4m)$ as $\sigma (t)$ and assume that such a $\sigma (t)$ satisfies Assumption \ref{A1}. 
\begin{Thm} \label{T2}
Let $j=2$ or $3$ and $(q,r)$ be an admissible pair with $n=j$. Then, for all  
$\phi_j \in L^{2}({\bf R}^j)$ and $F_j \in L_{\lambda}^{\tilde{q}',\tilde{r}'} $ with $n=j$, the following Strichartz estimates hold: 
\begin{align*}
\left\| 
U_{j,L}(t,0) \phi_j 
\right\|_{q,r, 2\lambda/j} \leq C \left\| \phi_j  \right\|_{2}, 
\end{align*}
\begin{align*}
\left\| 
\int_0^t U_{j,L}(t,0) U_{j,L}(s,0)^{\ast} F_j(s) ds 
\right\|_{q,r, 2\lambda/j } \leq C \left\| 
F_j
\right\|_{\tilde{q}', \tilde{r}', - 2\lambda/j}
\end{align*}
and 
\begin{align*}
\left\| 
\int_{-t}^0 U_{j,L}(t,0) U_{j,L}(s,0)^{\ast} F_j(s) ds 
\right\|_{q,r, 2\lambda/j} \leq C \left\| 
F_j
\right\|_{\tilde{q}', \tilde{r}', - 2\lambda/j}.
\end{align*}
\end{Thm}

Recently, nonlinear problems associated with \eqref{H0} have been considered by \cite{Ca} and Carles-Silva \cite{CaS}. In addition to these, Fanelli \cite{F} and Gonda-Machihara-Ozawa \cite{GMO} considered the nonlinear problems associated with the Hamiltonian $\hat{H}_0 (t) = -a(t) \Delta$. Using Korotyaev's transform, we can reduce the Hamiltonian \eqref{H0} to $\hat{H}_0 (t)$ with $a(t) = \zeta _1(t)^{-2}$ (see \cite{Ko} and \cite{KY}) and hence our result may be applicable to such studies. 

\section{Dispersive estimates}
To deduce the dispersive estimates for $ U_0(t,s)$, we use the following factorization, which was first obtained by \cite{Ko} and 
was rewritten by \cite{AK} through an{\em MDFM}-type decomposition (see (7.1) of \cite{AK});  
\begin{Lem}
Let $y_1 (t)$ and $y_2 (t)$ be the solutions to \eqref{2}. We define $a_1(t)$ and $a_2(t)$ as 
\begin{align} \label{3}
& \begin{cases}
a_1 (t) &= (y_1'(t) y_2 (t) -y_1 (t) y_2 '(t))/(y_1 (t)^2 + y_2 (t)^2), \\ 
a_1 (0) & = 1/m, 
\end{cases} \\ \label{4}
& \begin{cases}
a_2 (t) &= -(y_1'(t) y_1 (t) +y_2 (t) y_2 '(t))/(y_1 (t)^2 + y_2 (t)^2), \\ 
a_2 (0) & = 0. 
\end{cases}
\end{align}
Then, the following factorization of the propagator holds; 
\begin{align*}
U_0(t,0) = \CAL{M} \left( 
\frac{-1}{ma_2(t)}
\right) (i)^{n/2}\CAL{D} \left( 
\frac{1}{\sqrt{ma_1(t)}}
\right) e^{ -i \int_0^t a_1 (s) ds (-\Delta + x^2)/2},
\end{align*}
where 
$$ 
(\CAL{M}(\tau) \varphi ) (x) := e^{ix^2 /(2 \tau)} \varphi (x), \quad (\CAL{D}(\tau) \varphi) (x) := \frac{1}{(i \tau )^{n/2} } \varphi \left( \frac{x}{\tau} \right). 
$$
Moreover, the{\em MDFM}-type decomposition is 
\begin{align*}
U_0(t,0) &= \CAL{M} \left( 
\frac{-1}{ma_2(t)}
\right) (i)^{n/2}\CAL{D} \left( 
\frac{1}{\sqrt{ma_1(t)}}
\right) M\left( 
\tan \left( 
\int_0^t a_1(s) ds
\right)
\right) \\ & \quad \times D\left( 
\sin \left( 
\int_0^t a_1(s) ds
\right) 
\right) 
\SCR{F} M \left( 
\tan \left( 
\int_0^t a_1(s) ds
\right) 
\right)  , 
\end{align*}
where $\SCR{F}$ is the Fourier transform from $L^2({\bf R}^n)$ to $L^2({\bf R}^n)$.
\end{Lem}
Noting \eqref{2}, we have that for all $t \in {\bf R}$, 
\begin{align}\label{adad5}
 y_1'(t) y_2(t) -y_1(t) y_2 '(t) = W 
 \end{align} holds for some constant $W$. In this case, we notice $W= (y_1(0)^2 + y_2(0)^2) /m$ by \eqref{3}. 
Because we assumed that $y_1$ and $y_2$ are linearly independent, we have $W \neq 0$, which implies 
$a_1(t) = W /(y_1(t)^2+y_2(t) ^2 ) \geq 0$ holds. 

Here, we introduce estimations for coefficients that appear in the above factorization. 
\begin{Lem}
Let $a_1(t)$ and $a_2(t)$ be equivalent to those in \eqref{3} and \eqref{4}. We define $$ A(t) = \int_0^t a_1(s) ds. $$ Then, properties (i), (ii), and (iii) hold: \\ 
(i). For all sufficiently large constant $R>0$, there exists a constant $c_0 >0$ such that 
\begin{align} \label{7}
\sup_{|s| \leq R} (y_1(s)^2 +y_2(s)^2)^{-1} \leq c_0.
\end{align} 
(ii). For all sufficiently large $|t| \gg 1$, 
there exist constants
$0< c_{\mathrm{m}} < c_{\mathrm{M}}$ such that  
\begin{align} \label{6}
c_{\mathrm{m}} t^{2 \lambda -2} \leq a_1 (t) \leq c_{\mathrm{M}} t^{2 \lambda -2}.
\end{align}
\\ 
(iii). For all sufficiently large $|t| \gg 1$, there exists a constant $C_{\mathrm{A}}>0$ such that 
\begin{align} \label{m5}
A(t) \leq C_{\mathrm{A}}.
\end{align}
\end{Lem}
\Proof{ 
In this proof, we let $R \gg 1$ so that for all $|t| \geq R$, 
$2 |c_{1,\pm}||t|^{1- \lambda} \geq 
|y_1(t) |\geq |c_{1,\pm}||t|^{1- \lambda} /2 $ and $ 2|c_{2,\pm}||t|^{ \lambda} \geq |y_2 (t) |\geq |c_{2,\pm}||t|^{ \lambda} /2 $ hold and consider the case where 
$|t| \geq R$ and $s \leq R$.  

First, we prove (i). Roughly explaining, there exists $s_0 \in{\bf R}$ so that $(y_1(s_0)^2 +y_2(s_0))^{-1} = \infty$, which is equivalent to $y_1(s_0) =y_2(s_0) =0$. However, this contradicts \eqref{adad5} because $y_1'(t)$ and $y_2'(t)$ are continuous, i.e., are bounded for all $t \in[-R,R]$.

 Now, we prove \eqref{7} concretely. Suppose that for some $s_0 \in {\bf R}$, $|s_0| \leq R$, there exists a sufficiently small constant $\delta >0 $ such that 
$|y_1(s_0)| \leq \delta$. Then, noting $y_1(s_0) y_2'(s_0) - y_1'(s_0) y_2 (s_0) =W$, $\sup_{|s| \leq R} |y_1'(s)| \leq C_1$,
$\sup_{|s| \leq R} |y_2'(s)| \leq C_2$, and $\sup_{|s| \leq R} |y_1(s)| \leq C_3$, we have
\begin{align}\label{lo1}
 \delta C_2 \geq |y_1(s_0) y_2'(s_0)| = |W+ y_1'(s_0) y_2(s_0)| \geq W - W/2 \quad  & \quad \mbox{ if } |y_1'(s_0)| \ll 1, \\ 
 - \delta C_2 -W \leq y_1'(s_0) y_2(s_0) \leq -W +\delta C_2 \quad & \quad  \mbox{ if } |y_1'(s_0)| > 0 .\label{lo2}
\end{align}
By taking $\delta >0$ to be sufficiently small, inequality \eqref{lo1} fails because $W>0$, $C_2$, and $C_3$ are independent of $\delta$, which means $y_1'(s_0) \neq 0$ if $|y_1 (s_0)| \ll 1$. Moreover, by \eqref{lo2}, we obtain 
\begin{align*}
|y_2(s)| \geq (W-\delta C_2) /C_1.  
\end{align*}
This inequality means that if $y_1(s_0) =0$, then $|y_2(s)| \geq W/C_1$, i.e., $(y_1(s_0)^2 +y_2(s_0)^2)^{-1} 
\leq C_1^2W^{-2}$. For the case where 
$y_2(s_1) =0$, we can deduce $(y_1(s_1)^2 + y_2(s_1)^2)^{-1} \leq C_2 W^{-2}$ by almost the same calculation, which yields \eqref{7}. 

Next, we prove (ii). Using $W=y_1'(t) y_2(t) -y_1(t) y_2 '(t) > 0 $, 
$a_1 (t)$ can be written as 
\begin{align*}
 a_1 (t)  = (y_1(t))^{-2} \frac{W}{1 + (y_2(t)/y_1 (t))^2}.  
\end{align*}
This yields that for $|t| \geq R$ 
\begin{align*}
\frac{W}{2} |y_1(t)|^{-2} \leq \frac{W}{1+(4c_{2, \pm} t^{\lambda}/(c_{1, \pm} t^{1- \lambda}))^2} |y_1(t)|^{-2} \leq a_1(t) \leq  W |y_1(t)|^{-2}
\end{align*}
Combining \eqref{5} and this inequality, we have \eqref{6}. 

Finally, we prove (iii). Using \eqref{7}, we have 
\begin{align*}
A(t) &= W \left| \int_0^t \frac{ds}{y_1 (s)^2 + y_2 (s)^2 } \right|  
        \leq c_{{0}} RW + 4W \int_{R}^{\infty} s^{2\lambda -2} ds 
        \leq C_{\mathrm{A}},
\end{align*}
where we use $2 \lambda -2 < -1$ because $0 \leq \lambda <1/2$.
}

\begin{Lem} \label{L3}
For all $\phi \in L^1({\bf R}^n)$, 
\begin{align} \label{m2}
\left\| 
U(t,0)U_0(s,0)^{\ast} \phi 
\right\|_{L^{\infty}({\bf R}^n)} \leq C \frac{|a_1(t)a_1(s)|^{n/4}}{|\sin (A(t)-A(s))|^{n/2}} \left\| \phi \right\|_{L^1({\bf R}^n)}.
\end{align}
\end{Lem}
\Proof{
By the factorization of $U_0(t,0)$, for all $\phi \in \SCR{S}({\bf R}^n)$, 
\begin{align*}
& U_0(t,0)U_0(s,0)^{\ast} \phi \\ &= (-i)^{n/2} U_0(t,0) e^{-iA(s)(-\Delta +x^2)/2} \CAL{D}((ma_1(s))^{1/2})\CAL{M} ((ma_2(s))^{-1}) \phi \\ & 
=  (ma_1(s))^{-n/4} \CAL{M} (-(ma_2(t))^{-1}) \CAL{D}((ma_1(t))^{-1/2}) e^{-i(A(t)-A(s))(-\Delta +x^2)/2} \psi (s, (ma_2(s))^{-1/2} x)
\end{align*}
holds, where $\psi (s,x) = e^{ima_2(s)x^2/2} \phi(x)$. The integral kernel of $e^{-i\alpha(\Delta +x^2)/2}$ is well known (see, e.g.,  
(3.2) of \cite{AK}, H\"{o}rmander \cite{Ho}, and (2.10) of Bony-Carles-H\"{a}fner-Michel
 \cite{BCHM}) as 
\begin{align*}
\left| (e^{-i \alpha (-\Delta + x^2)/2}\varphi)(x) \right| &= \frac{1}{|2\pi i \sin \alpha|^{n/2}} \left|\int_{{\bf R}^n} e^{i (\cos \alpha (y^2 +x^2) - 2 x \cdot y)/(2 \sin \alpha)} 
\varphi (y) dy \right| \\ 
& =:  \frac{1}{|2\pi i \sin \alpha|^{n/2}} \left| \int_{{\bf R}^n} S(\alpha; x,y)
\varphi (y) dy \right|,
\end{align*}
which gives
\begin{align*}
\left| U_0(t,0) U_0(s,0)^{\ast} \phi \right| &=\left| e^{-ima_1(t)x^2/2} \frac{(ma_1(s) )^{-n/4} (ma_1(t))^{n/4}}{|2\pi\sin (A(t) - A(s))|^{n/2}} \right| \\ & \qquad \times \left| \int_{{\bf R}^n} 
S(A(t)-A(s); (ma_1(t))^{1/2}x, y )  \psi (s, (ma_2(s))^{-1/2} y) dy \right| .
\end{align*} 
By the change of a variable $(ma_2(s))^{-1/2} y = z$, we see that 
\begin{align*}
\left| U_0(t,0) U_0(s,0)^{\ast} \phi \right| &= \left| e^{-ima_1(t)x^2/2} \frac{(m^2 a_1(s) a_1(t))^{n/4}}{|2\pi \sin (A(t) - A(s))|^{n/2}} \right| \\ & \qquad \times \left|  \int_{{\bf R}^n} 
S(A(t)-A(s); (ma_1(t))^{1/2}x, (a_2(s)m)^{1/2}z )  \psi (s, z) dz \right|,
\end{align*}
which yields \eqref{m2}.

}

Next, we shall prove the dispersive estimate, which is equivalent to (1.10) and (1.11) of \cite{KY}. 

Let us define $r_0$ as a sufficiently large constant so that for all $|t| \geq r_0$, 
\begin{align*}
\frac{|c_{1, \pm}|}{2} |t|^{1- \lambda} \leq |y_1 (t)| \leq 2 |c_{1, \pm}| |t|^{1- \lambda}, \quad 
\frac{|c_{2, \pm}|}{2} |t|^{\lambda} \leq  |y_2 (t)| \leq 2 |c_{2, \pm}| |t|^{ \lambda}
\end{align*} 
and 
\begin{align}\label{8}
\left| \int_{\pm r_0}^{\pm \infty}  \frac{W}{y_1(\tau) ^2 + y_2 (\tau) ^2} d \tau \right| < \pi /2
\end{align}
hold. Then, for such a $r_0$, we define the spaces $\Omega_{0, \pm}$ and $\Omega _{\lambda, \pm}$ as follows: 
\begin{align*}
\Omega _{0, \pm} := \left\{ 
(t,s) \in {\bf R} \times {\bf R} \backslash [-r_0,r_0] ^2 \, \Big| \, |t|/2 \leq |s| \leq 2 |t|, \ \pm t \geq r_0, \, \pm s \geq r_0
\right\}
\end{align*}
and 
\begin{align*} 
\Omega _{\lambda, \pm} := \left\{ 
(t,s) \in {\bf R} \times {\bf R} \backslash [-r_0,r_0] ^2 \, \Big| \, |s| \leq |t|/2 \ \mbox{or} \  |s| \geq 2 |t|, \ \pm t \geq r_0, \, \pm s \geq r_0
\right\}.
\end{align*}

\begin{Lem} \label{L1}
On $\Omega_{0,\pm}$, 
\begin{align} \label{11}
\left\| U_0(t,0) U_0(s,0) ^{\ast} \phi \right\|_{L^q({\bf R}^n)} \leq 
C |t-s|^{-n(1/2-1/q)} \left\| \phi \right\|_{L^{q'}({\bf R}^n)}
\end{align} 
holds. Moreover, on $\Omega_{\lambda , \pm} $, 
\begin{align} \label{12}
\left\| U_0(t,0) U_0(s,0) ^{\ast} \phi \right\|_{L^q({\bf R}^n)} \leq 
C |t-s|^{-n(1- \lambda)(1/2-1/q)} \left\| \phi \right\|_{L^{q'}({\bf R}^n)}
\end{align} 
holds.
\end{Lem}
\Proof{
Before proving this lemma, note that for $|t| \geq r_0$, we have \eqref{6}. 
First, consider the case where $(t,s) \in \Omega _{0, \pm}$ and $|t-s| \leq 1$. Then, by the mean value theorem, there exists $\tau \in [s,t]$ or $[t,s]$ such that 
\begin{align*} 
& |t-s|^{n/2} |a_1(t) a_1 (s)|^{n/4} |\sin (A(t) -A(s))|^{-n/2} \\ & \leq \left| 
\frac{A(t) - A(s)}{\sin (A(t) -A(s))}
\right|^{n/2} \left| 
A'(\tau)
\right| ^{-n/2}   |a_1(t) a_1 (s)|^{n/4}
\end{align*}
holds. Because 
\begin{align*}
A'(\tau) = a_1(\tau) \geq c_{\mathrm{m}} |\tau|^{-2 +2 \lambda}\geq 
\begin{cases} 
 c_{\mathrm{m}} |t|^{-2 +2 \lambda} & |t| < |s| , \\
 c_{\mathrm{m}} |s|^{-2 +2 \lambda} & |s| < |t| 
 \end{cases} 
\end{align*}
and \eqref{6}, we have 
\begin{align} \label{9}
|t-s|^{n/2} |a_1(t) a_1 (s)|^{n/4} |\sin (A(t) -A(s))|^{-n/2} \leq C.
\end{align}
For $(t,s) \in \Omega _{0, \pm}$, $t>s$, and $|t-s| > 1$, noting \eqref{8}, we have $0< A(t) -A(s) < \pi/2$, which gives 
\begin{align} \nn
|t-s|^{n/2} |a_1(t) a_1 (s)|^{n/4} |\sin (A(t) -A(s))|^{-n/2} &\leq C(|t|^{n/2} +|s|^{n/2}) |t|^{-n(1- \lambda)} \\ &\leq C|t|^{-(n/2)(1-2 \lambda)} \leq C. 
\label{10}
\end{align}
Inequalities \eqref{9} and \eqref{10} mean that
\begin{align*}
|a_1(t) a_1(s)|^{n/4} |\sin (A(t)-A(s))| ^{-n/2} \leq C|t-s|^{-n/2}
\end{align*}
holds for all $(t,s) \in \Omega_{0,\pm}$.

Next, we consider the case where $(t,s) \in \Omega _{\lambda, \pm}$ and $|s| \leq |t|/2$. Here, we remark that on $(t,s) \in \Omega _{\lambda , \pm}$, it 
follows that there exists $C>0$ so that $|\sin (A(t)-A(s))| \geq C$ holds by \eqref{8}. Then, using $|t| \geq |t-s +s| \geq |t-s|-|t|/2$, i.e., 
$|t| \geq 2|t-s| /3$, we have 
\begin{align*}
|a_1(t) a_1 (s)|^{n/4} |\sin (A(t) -A(s))|^{-n/2} \leq C |t|^{-n(1-\lambda)/2} \leq C|t-s|^{-n(1- \lambda)/2}. 
\end{align*}
In the same way, on $(t,s) \in \Omega _{\lambda, \pm}$ and $2|t| \leq |s|$, we have 
$$ 
|a_1(t) a_1 (s)|^{n/4} |\sin (A(t) -A(s))|^{-n/2} \leq C|t-s|^{-n(1- \lambda)/2}. 
$$
These inequalities yield 
\begin{align*}
& \left\| U_0(t,0) U_0(s,0) ^{\ast} \phi \right\|_{L^{\infty}({\bf R}^n)} \leq 
C |t-s|^{-n/2} \left\| \phi \right\|_{L^{1}({\bf R}^n)}, \quad \mbox{on} \quad \Omega_{0,\pm}, \\
& \left\| U_0(t,0) U_0(s,0) ^{\ast} \phi \right\|_{L^{\infty}({\bf R}^n)} \leq 
C |t-s|^{-n(1- \lambda)/2} \left\| \phi \right\|_{L^{1}({\bf R}^n)}, \quad \mbox{on} \quad \Omega_{\lambda , \pm}, 
\end{align*} 
hold. These inequalities, $\left\|  U_0(t,0) U_0(s,0) ^{\ast} \phi  \right\|_{L^2({\bf R}^n)} =\left\| \phi  \right\|_{L^2({\bf R}^n)}  $, and 
Riesz-Thorin's interpolation theorem yield \eqref{11} and \eqref{12}.

}

\section{Proof of Theorem \ref{T1}}
In this section, we will prove Theorem \ref{T1}. The fundamental approach of the proof is based on \cite{KY}; see also Them 2.3.3 in Cazenave \cite{C}. 
To prove the Strichartz estimates 
for our model (except for the end-point),  
it is enough to prove the following inequalities:  
\begin{align} \label{13}
\left\| 
\int_{-\infty}^{\infty} U(t,0)^{\ast} F(t) dt 
\right\|_{2} \leq C \left\| F \right\|_{q',r',- \lambda}
\end{align} 
and 
\begin{align} \label{14}
\left\| 
\int_0^{|t|} U_0(t,0) U_0(s,0)^{\ast} F(s) ds 
\right\| _{q,r,\lambda} \leq C \left\| F \right\| _{q',r', -\lambda}
\end{align}
for an admissible pair $(q,r)$ and $F \in L_{-\lambda}^{q',r'}$ (see \cite{KY}), where we remark that $\left\| \cdot \right\|_{q',r', \lambda} \leq 
\left\| \cdot \right\|_{q',r', 0}$. However, by the singularity of 
$( \sin (A(t) -A(s) ))^{-1}$ (res. the singularity of $(\zeta _1(t) \zeta _2(s) - \zeta _1 (s) \zeta _2 (t))^{-1}$ in \cite{KY}), 
$L_{-\lambda}^{q',r'} $ needs to be divided into 
$L_{-\lambda, 0, \pm}^{q',r'} + L_{-\lambda, \ast}^{q',r'} $ with 
\begin{align*}
L_{-\lambda, 0, \pm}^{q',r'}  &:= \left\{ 
F \in L_{-\lambda}^{q',r'} \, \Big| \, F(s) = 0 \mbox{ for } \pm s \leq r_0 
\right\}, \\ 
L_{-\lambda, \ast }^{q',r'}  &:= \left\{ 
F \in L_{-\lambda}^{q',r'} \, \Big| \, F(s) = 0 \mbox{ for }  s \leq - r_0 \wedge r_0 \leq s
\right\}, 
\end{align*}  
to prove \eqref{13} and \eqref{14}.

By Lemma \ref{L1} and a completely identical approach to that in \cite{KY}, we can prove \eqref{13} and \eqref{14} for $F \in L_{-\lambda,0,\pm}^{q',r'}$; hence, 
in the following, we only consider the case of $F \in L_{- \lambda, \ast}^{q',r'}$. 

First, we consider the case $(t,s) \in [-r_0,r_0] ^2 
$. Here, we remark that 
\begin{align*}
\sup_{\tau \in [-r_0,r_0]} a_1(\tau) \leq c_0 W
\end{align*}
holds. Let us define $\tilde{N} \in {\bf N}$ and suppose $r_0 >0$, such that 
\begin{align}\label{ax1}
\tilde{N}\pi < \int_{-\infty}^{\infty} \frac{W}{y_1 ^2(\tau ) + y_2^2 (\tau)}d\tau \leq  (\tilde{N} +1) \pi, \quad \int_{-\infty}^{r_0} \frac{W}{y_1 ^2(\tau ) + y_2^2 (\tau)}d\tau > \tilde{N} \pi 
\end{align} 
where we note, by \eqref{m5}, $\tilde{N} $ is a finite number. For some $\tilde{\delta} >0$ with $\tilde{\delta} < \pi /4$, suppose that $r_0 >0$ satisfies 
\begin{align}\label{adad10}
\int_{r_0}^{\infty}  \frac{W}{y_1 ^2(\tau ) + y_2^2 (\tau)}d\tau \geq 2 \tilde{\delta}, \quad 
\int_s^{\infty} \frac{W}{y_1 ^2(\tau ) + y_2^2 (\tau)}d\tau < (\tilde{N} +1) \pi
 .
\end{align} 
holds for all $s \in [-r_0,r_0]$. We define $N \in {\bf N} \cup \{0\}$ with $N \leq \tilde{N} $, and $\Omega _{\mathrm{res},N}$ and $\Omega _{\mathrm{reg}}$ as follows 
\begin{align*}
\Omega _{\mathrm{res},N} := \left\{ 
(t,s) \in [-r_0,r_0]^2 \, ; \, N \pi - \tilde{\delta} \leq |A(t) -A(s) |  \leq N \pi +\tilde{\delta} 
\right\}
\end{align*}
and 
\begin{align*}
\Omega _{\mathrm{reg}} := [-r_0,r_0] ^2 \backslash \bigcup_{N=0}^{\tilde{N}} \Omega _{\mathrm{res}, N}.
\end{align*}

Now, we investigate $\Omega _{\mathrm{res},N}$ more precisely. We only consider the case where $N \geq 1$ because the case of $N=0$ is easier to consider. \\ 
{\bf I}. For the case where $t>s$. \\ 
We define $\omega _N (s) $ with $\omega _N \geq 0$ and $s^{\ast} \in [-r_0,r_0]$ as follows 
\begin{align} \label{16}
\int_s^{s + \omega _N (s)}a_1(\tau) d \tau = N \pi, \quad 
\int_{s^{\ast}}^{\infty} a_1(\tau) d \tau = N \pi  , 
\end{align}
i.e., $s^{\ast}+ \omega_N(s^{\ast}) = \infty$ (blow-up point). 
Below, we shall prove that $[-r_0,r_0] \times [s^{\ast}, r_0] \not\subset \Omega_{\mathrm{res},N}$ and there exists a constant $\tilde{c}_0 >0$ such that 
\begin{align} \nn
|\sin (A(t) -A(s))| &=\left| \frac{\sin (A(t) -A(s) -N\pi)}{A(t)-A(s)-N \pi} \frac{A(t)-A(s) -N \pi}{t-s-\omega _N(s)} (t-s-\omega _N(s))
\right| \\ & \geq 
\tilde{c}_0  |t-s-\omega _N (s)| \label{mk6}
\end{align}
holds on $(t,s) \in \Omega_{\mathrm{res}, N} \cap \{ t>s\}$. 

Here, we remark that $t \leq r_0$. Then, for $s^{\ast} \leq s \leq r_0$, we have   
\begin{align*}
A(t)-A(s)
 = 
\int_s^t a_1(\tau) d \tau
 \leq 
 \int_{s^{\ast}}^{r_0} a_1(\tau) d \tau 
 \leq  \int_{s^{\ast} }^{\infty} a_1(\tau) d \tau - \int_{r_0}^{\infty}a_1(\tau) d\tau \leq N \pi-2 \tilde{\delta}
\end{align*} 
which holds by \eqref{adad10} and \eqref{16}. By the definition of $\Omega_{\mathrm{res},N}$, we have that 
\begin{align*}
\left\{ 
(t,s) \in [-r_0,r_0]^2 \, | \,  t > s , \ s\geq s^{\ast} \right\} 
 \not\subset \Omega_{\mathrm{res}, N}
\end{align*}
holds. Hence we have that $\Omega_{\mathrm{res},N}\cap \{ t>s \} \subset [-r_0,r_0] \times [-r_0,s^{\ast \ast}] $ with 
\begin{align*}
 s^{\ast \ast} := \max_{s \in [-r_0,s^{\ast})} \left\{ \int_{s^{}}^{\omega_N(s^{}) + s^{}} a_1(\tau) = N \pi- \tilde{\delta} \right\}, \quad \omega_N(s^{\ast\ast}) + s^{\ast \ast} < \infty.
\end{align*}
and obtain 
\begin{align*}
\inf_{\tau \in [-r_0,  \omega_N(s^{\ast \ast})+s^{\ast \ast} ]} a_1(\tau) \geq c_0 
\end{align*} 
for some $c_0 >0$. 
By simple calculation, we have  
\begin{align} 
 A(t) -A(s) -N \pi &= \int_s^t  a_1(\tau ) d\tau - \int_s^{s + \omega _N (s)} a_1(\tau) d \tau 
 = \int_{s + \omega _N (s)}^t a_1(\tau) d \tau . \label{15}
\end{align} and we get 
\begin{align*}
|A(t) -A(s) -N \pi | \geq \int_{s+\omega_N(s)}^t a_1(\tau) d \tau \geq C(s^{\ast\ast}) |t-s- \omega_N(s)| . 
\end{align*} 
These show that \eqref{mk6} holds on $(t,s) \in \Omega _{\mathrm{res},N} \cap \{t>s \}$.
 \\
{\bf II}. For the case where $t \leq s $. \\ 
By the same calculation in \eqref{15}, we have 
\begin{align*}
A(t)-A(s) + N \pi = \int_s^t a_1 (\tau) d \tau + \int^{t + \omega _N(t) }_t a_1(\tau) d \tau = \int_{s}^{t+ \omega_N(t)} a_1(\tau) d \tau .
\end{align*}
For the same reason, we get 
\begin{align*}
\left\{ 
(t,s) \in [-r_0,r_0] ^2 \, | \, t \leq s , \ t \geq t^{\ast}
\right\} \not\subset \Omega_{\mathrm{res},N}, \quad \int_{t^{\ast}}^{\infty} a_1(\tau) d \tau = N \pi. 
\end{align*} 
Hence, we have 
\begin{align*}
|A(t) -A(s) + N \pi | \geq \int^{t+\omega_N(t)}_s a_1(\tau) d \tau \geq C(t^{\ast\ast}) |t+\omega_N(t)-s|,  
\end{align*} 
where 
\begin{align*}
 t^{\ast \ast} := \max_{t \in [-r_0,t^{\ast})} \left\{ \int_{t^{}}^{t+ \omega_N(t^{})} a_1(\tau) = N \pi- \tilde{\delta} \right\}, \quad \omega_N(t^{\ast\ast}) + t^{\ast \ast} < \infty.
\end{align*}
These show that there exists a constant $\tilde{c}_0 >0$ such that \eqref{mk6} holds on $(t,s) \in \Omega_{\mathrm{res}, N} \cap \{ t\leq s\}$.

Now, we estimate $\left\| U_0(t,0)^{\ast} F(t) \right\|_{2}^2$ for $F \in
L_{-\lambda, \ast}^{q',r'}$. Let us define 
$$ 
\theta _{\lambda} := (1-\lambda) \left( \frac12 - \frac{1}{r} \right)
$$ 
for some $r \geq 2$. By \eqref{m2} and straightforward calculation, 
\begin{align*}
& \left\| U_0(t,0)^{\ast} F(t) \right\|_{2}^2  \\ & \leq C \iint _{(t,s)
 \in [-r_0,r_0]^2}
 |a_1(t)a_1(s)|^{n\theta _0/2} |\sin (A(t) -A(s))|^{-n\theta _{0}} \left\| 
F(s)
\right\|_{r'} \left\| 
F(t)
\right\|_{r'} ds dt
\\ & \leq C \iint_{\Omega_{\mathrm{reg}} } |\sin (A(t)-A(s))| ^{-n\theta _0} \left\| 
F(s)
\right\|_{r'} \left\| 
F(t)
\right\|_{r'} dsdt
\\ & \quad + C \sum_{N=0}^{\tilde{N}} \iint_{\Omega_{\mathrm{res},N}} |\sin (A(t)-A(s))| ^{-n\theta _0} \left\| 
F(s)
\right\|_{r'} \left\| 
F(t)
\right\|_{r'} dsdt \\ &=:
I_{\mathrm{reg}} + I_{\mathrm{res}} .
\end{align*} 
holds. Because $|\sin (A(t)-A(s))| \geq C$ and $|t-s| \leq 2r_0$ hold on $\Omega _{\mathrm{reg}}$, we have 
\begin{align}
I_{\mathrm{reg}} \leq  C \iint_{\Omega_{\mathrm{reg}}} |t-s| ^{-n\theta _0} \left\| 
F(s)
\right\|_{r'} \left\| 
F(t)
\right\|_{r'} dsdt \leq C \left\| 
F
\right\|_{q',r',0}, \label{m10}
\end{align}
where we use the Hardy-Littlewood-Soblev inequality. 

Now, we estimate $I_{\mathrm{res}}$. For 
$$\Omega_{>,N}:=
\{(t,s)\in \Omega_{\mathrm{res},N} \, | \, t>s  \} \subset \{(t,s)\in[-r_0,r_0] \times [-r_0, s^{\ast \ast}]  \, | \, t>s  \}$$ and $$\Omega_{\leq,N} :=
\{(t,s)\in \Omega_{\mathrm{res},N} \, | \, t \leq s \} \subset
\{(t,s)\in [-r_0, t^{\ast \ast}] \times [-r_0,r_0] \, | \, t \leq s \} , $$ we devide
$I_{\mathrm{res}}$ into $I_{\mathrm{res},>} + I_{\mathrm{res},\leq}$
with 
\begin{align*}
I_{\mathrm{res},>} &:=  C \sum_{N=0}^{\tilde{N}} \iint_{\Omega_{>,N}} |\sin (A(t)-A(s))| ^{-n\theta _0} \left\| 
F(s)
\right\|_{r'} \left\| 
F(t)
\right\|_{r'} dsdt, \\ 
I_{\mathrm{res},\leq} &:= C \sum_{N=0}^{\tilde{N}} \iint_{\Omega_{\leq,N}} |\sin (A(t)-A(s))| ^{-n\theta _0} \left\| 
F(s)
\right\|_{r'} \left\| 
F(t)
\right\|_{r'} dsdt. 
\end{align*} 
We only estimate $I_{\mathrm{res},>}$ because $I_{\mathrm{res},\leq} $ can
be estimated in the same way. By \eqref{mk6} and \eqref{15}, we have 
\begin{align}
 I_{\mathrm{res},>} \leq  C 
 \sum_{N=0}^{\tilde{N}} \iint_{\Omega_{>,N}} |t-s-\omega _N(s)| ^{-n\theta _0} \left\| 
F(s)
\right\|_{r'} \left\| 
F(t)
\right\|_{r'} dsdt .\label{m8}
\end{align} 
Because \eqref{16}, we have  
\begin{align*}
0 = \frac{d}{ds} \int_s^{s + \omega _N(s)}a_1(\tau) d \tau =(1+\omega _N'(s)) a_1(s+\omega _N(s)) -a_1(s), 
\end{align*}
and that shows that there exists $\tilde{c} >0$ such that
\begin{align}
1+ \omega_N '(s) = \frac{a_1(s)}{a_1(s+ \omega _N(s))}  \geq \frac{\tilde{c}}{c_0 W}, \label{m7}
\end{align}
i.e., $f(s) = s+\omega_N(s)$ is a monotonic increasing function. Noting $\Omega _{>,N} \subset [-r_0,r_0] \times [-r_0, s^{\ast\ast} ]$, \eqref{m8}, \eqref{m7} and a change of variables 
\begin{align*} 
\tilde{f} (s) := s+ \omega _N(s) = \tau,  \quad \mbox{for} \   -r_0 \leq s \leq s^{\ast\ast}, \end{align*}
 we also have 
\begin{align*}
 I_{\mathrm{res},>}  & \leq C \left( \sup_{-r_0 \leq s \leq s^{\ast\ast} } |1+\omega _N'(s)|^{-1} \right) \sum_{N=0}^{\tilde{N}} \int_{-r_0}^{r_0  } \int_{-r_0+\omega_N(-r_0)}^{s^{\ast\ast}+\omega _N({s}^{\ast \ast})} \\ & \qquad \times\Big\{ |t-\tau| ^{-n\theta _0} \left\| 
F(\tilde{f}^{-1}{(\tau)})
\right\|_{r'} \left\| 
F(t)
\right\|_{r'} \Big\} d\tau dt \\ 
& \leq C \left\| 
F
\right\|_{q',r',0}, 
\end{align*} 
where $\tilde{f}$ is a bijection from $[-r_0,s^{\ast\ast}]$ to
$[f(-r_0), f(s^{\ast\ast})]$. Hence, we have \eqref{13} for $F \in L_{- \lambda , \ast}^{q',r'}$. 

Next, we shall prove \eqref{14} for $F \in L_{- \lambda , \ast, \pm}^{q',r'}$. 
We define $\tilde{\Omega}_{\mathrm{res},N} (t)$ and $\tilde{\Omega} _{\mathrm{reg}} (t)$ as follows: 
\begin{align*}
\tilde{\Omega} _{\mathrm{res},N}(t) := \left\{ 
s \in [-r_0,r_0] \, ; \, N \pi - \tilde{\delta} \leq |A(t) -A(s) | \leq N \pi +  \tilde{\delta}, 
\right\}
\end{align*}
and 
\begin{align*}
\tilde{\Omega} _{\mathrm{reg}} (t) :=  [-r_0,r_0]  \backslash \bigcup _{N=0}^{\tilde{N}} \tilde{\Omega}  _{\mathrm{res}, N}(t)  .
\end{align*}
Here, we remark that on $t >r_0$ and $s \in \tilde{\Omega}_{\mathrm{reg}}(t)$, 
\begin{align*}
& |a_1(t)a_1(s)|^{1/2} |\sin (A(t)-A(s))|^{-1} \leq C |a_1(t)|^{1/2} \leq C |t|^{-(1- \lambda)} \\ &= 
C |t-s|^{1- \lambda}  |t|^{-(1- \lambda)} |t-s|^{-(1- \lambda)} \leq C |t-s|^{-(1- \lambda)}
\end{align*}
holds; on $0 \leq t \leq r_0$ and $s \in \tilde{\Omega}_{\mathrm{reg}}(t)$, 
\begin{align*}
& |a_1(t)a_1(s)|^{1/2} |\sin (A(t)-A(s))|^{-1} \leq C  \leq C |t-s|^{-(1- \lambda)}  
\end{align*}
holds, where we use $|t-s| \leq 2r_0$. Then, using the same scheme as in the proof of \eqref{14}, we have 
\begin{align*}
& \left\| 
\int_0^{|t|} U_0(t,0) U_0(s,0)^{\ast} F(s) ds 
\right\| _{q,r,\lambda}^q \\  &\leq C\int_{-\infty}^{\infty} \left( \int_{-r_0}^{r_0}
|a_1(t) a_1(s)|^{n\theta _0/2} |\sin (A(t)-A(s))|^{-n \theta _{0}}
 \left\| F(s) \right\| _{r'} ds \right) ^q dt \\ & \leq 
C \int_{-\infty}^{\infty} \left( \int_{s\in \tilde{\Omega}_{\mathrm{reg}} (t)}
|t-s|^{-n \theta _{\lambda}}
 \left\| F(s) \right\| _{r'} ds \right) ^q dt \\ & \quad + 
 C\sum_{N=0}^{\tilde{N}} \int_{-\infty}^{\infty} \left( \int_{s\in \tilde{\Omega}_{\mathrm{res},N}(t) } 
|a_1(t) a_1(s)|^{n\theta _0/2} |\sin (A(t)-A(s))|^{-n \theta _{0}}
 \left\| F(s) \right\| _{r'} ds \right) ^q dt \\ 
 & =: J_{\mathrm{reg}} + J_{\mathrm{res}}. 
\end{align*}
Clearly, 
\begin{align*}
& J_{\mathrm{reg}}  \leq C \int_{-\infty}^{\infty} \left( \int_{s\in \tilde{\Omega}_{\mathrm{reg}} (t)}
|t-s|^{-n \theta _{\lambda}}
 \left\| F(s) \right\| _{r'} ds \right) ^q dt \leq C \| F \|^q_{q',r',0} \leq C \| F \|^q _{q',r',-\lambda}
\end{align*}
holds by the Hardy-Littlewood-Soblev inequality. Now, we estimate $J_{\mathrm{res}}$. First, we divide the limits of integration into 
\begin{align*}
J_{\mathrm{res}} = C \sum_{N=0}^{\tilde{N}} \left(  J_{\mathrm{res}, \mathrm{in
},N} +  J_{\mathrm{res},\mathrm{out},N} \right) 
\end{align*} 
with 
\begin{align*} 
J_{\mathrm{res}, \mathrm{in} ,N} := \int_{|t|  \leq {r}_0} \left( \int_{s\in \tilde{\Omega}_{\mathrm{res},N}(t) } 
|a_1(t) a_1(s)|^{n\theta _0/2} |\sin (A(t)-A(s))|^{-n \theta _{0}}
 \left\| F(s) \right\| _{r'} ds \right) ^q dt
\end{align*} 
and 
\begin{align*} 
J_{\mathrm{res}, \mathrm{out} ,N} := \int_{|t|  \geq {r}_0} \left( \int_{s\in \tilde{\Omega}_{\mathrm{res},N}(t) } 
|a_1(t) a_1(s)|^{n\theta _0/2} |\sin (A(t)-A(s))|^{-n \theta _{0}}
 \left\| F(s) \right\| _{r'} ds \right) ^q dt. 
 \end{align*}
 Then, with the same calculation as in the proof of \eqref{13}, we obtain 
 \begin{align*}
C \sum_{N=0}^{\tilde{N}}  J_{\mathrm{res}, \mathrm{in} ,N}  \leq C \| F \|^q _{q',r',-\lambda} .
 \end{align*}  
Now, we estimate $J_{\mathrm{res}, \mathrm{ out} ,N}$. \\ 
For $t \geq r_0$. We define $\tilde{\omega} _N (t) \geq 0$ so that 
\begin{align*}
\int^{t}_{t-\tilde{\omega} _N(t)}  a_1(\tau) d \tau = N \pi.   
\end{align*} 
By the assumption \eqref{ax1}, we have that 
\begin{align*}
\int_{-\infty}^{t} a_1(\tau) d \tau > N \pi, \mbox{ for } t>r_0
\end{align*} 
holds for all $N \leq \tilde{N}$, which means that on $t \geq {r}_0$, $\tilde{\omega}_N(t)$ has no blow-up points ($\omega _N(t_0) = - \infty$ for some $t_0 \in (- \infty, \infty)$). Then, with the same calculation as in \eqref{15} and \eqref{mk6}, for the case of $t \geq {r}_0 $, there exists $\tilde{c}_2 >0$ such that 
\begin{align*}
\left| 
\sin (A(t)-A(s))
\right| \geq \tilde{c}_2 |t-\tilde{\omega} _N (t) -s|. 
\end{align*}
Then, it holds that, for some compact set $\hat{C}$,  
\begin{align*}
& \int_{t  \geq {r}_0} \left( \int_{s\in \tilde{\Omega}_{\mathrm{res},N}(t)}  
|a_1(t) a_1(s)|^{n\theta _0/2} |\sin (A(t)-A(s))|^{-n \theta _{0}}
 \left\| F(s) \right\| _{r'} ds \right) ^q dt \\ & 
\leq  C \int_{t  \geq {r}_0} \left( \int_{s\in \tilde{\Omega}_{\mathrm{res},N}(t)}  
|a_1(t) |^{n\theta _0/2} |t-\tilde{\omega} _N (t) -s|^{-n \theta _{0}}
 \left\| F(s) \right\| _{r'} ds \right) ^q dt 
 \\ & 
 \leq 
 C \int_{\tau \in \hat{C}} 
 \sup_{t \geq {r}_0} \left\{ |a_1(t)|^{qn\theta _0/2 -1 } |a_1(t-\omega_N(t))|^{} \right\}\left( 
 \int_{{-r}_0}^{r_0} \left| \tau -s \right|^{-n \theta _0} \left\| F(s) \right\| _{r'} ds 
 \right) ^q d\tau
\end{align*} 
 holds, where we use $1-\tilde{\omega}_N'(t) = a_1(t) (a_1(t- \tilde{\omega}_N(t)))^{-1}$. By 
 \begin{align*}
 \frac{qn\theta _0}{2} -1 = q \left( 
 \frac{n}{4} - \frac{n}{2r}
 \right) -1 = 0, 
 \end{align*}
the above inequality is smaller than 
 \begin{align*}
 C \int_{\tau \in \hat{C}} 
 \left( 
 \int_{{-r}_0}^{r_0} \left| \tau -s \right|^{-n \theta _0} \left\| F(s) \right\| _{r'} ds 
 \right) ^q d\tau \leq C \| F \|_{q',r',0}^q \leq C \| F \|_{q',r',- \lambda}^q 
 \end{align*}
 For $t \leq -r_0$, we again employ $\omega_N(t)$ and consider 
 \begin{align*}
 A(t) -A(s) + N \pi = \int_s^{t+\omega_N(t)} a_1(\tau) d \tau. 
 \end{align*}
 Then, for some compact set $\hat{\hat{C}}$, we have 
\begin{align*}
& \int_{t\leq -r_0} \left( \int_{s\in \tilde{\Omega}_{\mathrm{res},N}(t)}  
|a_1(t) a_1(s)|^{n\theta _0/2} |\sin (A(t)-A(s))|^{-n \theta _{0}}
 \left\| F(s) \right\| _{r'} ds \right) ^q dt \\ & \leq 
C \int_{\tau \in \hat{\hat{C}}} 
 \left( 
 \int_{{-r}_0}^{r_0} \left| \tau -s \right|^{-n \theta _0} \left\| F(s) \right\| _{r'} ds 
 \right) ^q d\tau \leq C \| F \|_{q',r',0}^q \leq C \| F \|_{q',r',- \lambda}^q .
\end{align*}

Consequently, we have \eqref{13} and \eqref{14} for all $F \in L_{- \lambda}^{q',r'}$. These results and the same approach as that used in \cite{KY} 
yield Theorem \ref{T1}.

\section{Time-Decaying Magnetic Fields}
In this section, we shall consider the case of time-decaying magnetic fields. The Hamiltonian for this case can be written as in \eqref{Mg01} or \eqref{Mg02}. We define $U_{j,L} (t,0)$ as a propagator for $H_{j,L}(t)$, $j \in \{ 2,3\}$ and suppose $\sigma (t) = q^2B(t)^2/(4m)$ satisfies Assumption \ref{A1}. Summarizing these, we obtain for $j \in \{ 2,3 \}$
\begin{align}\label{Mg03}
H_{j,L} = - \frac{\Delta _{(2)}}{2m} - (j-2) \frac{\partial _{3}^2}{2m} + \frac{q^2B(t)^2}{8m} |x_{(2)}|^2 -\frac{qB(t)}{2m} L_{(2)} , \quad \mbox{ on } L^2({\bf R}^j).
\end{align} 
Because the operator $L_{(2)}$ commutes with other operators in \eqref{Mg03}, the propagator $U_{j,L} (t,0)$ can be reduced to 
\begin{align*}
U_{j,L} (t,0) = e^{i\Omega (t) L_{(2)}} \tilde{U}_{0} (t,0), 
\end{align*}
where $\Omega (t) = {\displaystyle \int_0^t qB(s)/(2m) ds}$ and $\tilde{U}_0(t,0)$ is a propagator for 
\begin{align*}
\tilde{H}_0(t) :=  - \frac{\Delta _{(2)}}{2m} - (j-2) \frac{\partial _{3}^2}{2m} + \frac{q^2B(t)^2}{8m} |x_{(2)}|^2.
\end{align*}
Hence, using the same calculations as in Lemma \ref{L1}, the following Lemma can be proven; this lemma immediately provides the proof of Theorem \ref{T2}. 
\begin{Lem}
Let $j = 2$ or $3$, $\phi_j ({\bf R}^j) \in \SCR{S}({\bf R}^j)$. On $\Omega_{0, \pm}$, the following estimate holds 
\begin{align*}
\left\| 
U_{j, L} (t,0)U_{j,L} (s,0)^{\ast} \phi_j 
\right\|_{q} \leq C |t-s|^{-j(1/2-1/q)} \left\| \phi_j \right\|_{q'} 
\end{align*}
and on $\Omega_{\lambda, \pm}$, the following estimate holds 
\begin{align*}
\left\| 
U_{j, L} (t,0)U_{j,L} (s,0)^{\ast} \phi_j 
\right\|_{q} \leq C |t-s|^{-j(1-(2 \lambda)/j)(1/2-1/q)} \left\| \phi_j \right\|_{q'}.
\end{align*}
\end{Lem}
\proof{
We define $\CAL{M}_{(2)}$ and $\CAL{D}_{(2)}$ as 
\begin{align*} 
& (\CAL{M}_{(2)}(\tau) \varphi ) (x_1,x_2,x_3) := e^{i(x_1^2+x^2_2) /(2 \tau)} \varphi (x_1,x_2,x_3), \\
& (\CAL{D}_{(2)}(\tau) \varphi) (x_1,x_2,x_3) := \frac{1}{i \tau } \varphi \left( \frac{x_1}{\tau}, \frac{x_2}{\tau},x_3 \right).
\end{align*}
Because $L_{(2)}$ commutes with $\CAL{M}_{(2)}(\cdot)$, $\CAL{D}_{(2)}(\cdot)$, and $-\Delta_{(2)} + x_1^2 +x_2^2 + (- \partial _3^2)$, we have 
\begin{align*}
& U_{j,L} (t,0) \\ &= e^{i(j-2) t \partial _3^2 /(2m)} e^{i\Omega (t) L_{(2)}} \CAL{M}_{(2)} \left( 
\frac{-1}{ma_2(t)}
\right) i\CAL{D}_{(2)} \left( 
\frac{1}{\sqrt{ma_1(t)}}
\right) e^{-i \int_0^t a_1 (s) ds (-\Delta _{(2)} + x_1^2+x_2^2)/2}, 
\end{align*}
Hence, using the same calculation as in the proof of Lemma \ref{L3}, we have 
\begin{align*}
& U_{3,L} (t,0) U_{3,L}(s,0)^{\ast} \phi_3 
= e^{-ima_1(t)(x_1^2+x_2^2) /2} \frac{(m^2 a_1(s) a_1(t))^{1/2}}{2\pi \sin (A(t) - A(s))} e^{i (t-s) \partial _3^2 /(2m)} e^{i\Omega (t-s) L_{(2)}}
\\ & \times\int_{{\bf R}^n} 
S(A(t)-A(s); (ma_1(t))^{1/2}x_1,ma_1(t))^{1/2}x_2 ; (a_2(s)m)^{1/2}z_1, (a_2(s)m)^{1/2} z_2 ) \\ & \qquad \times \psi_3 (s, z_1,z_2,x_3) dz,
\end{align*}
where $\psi_3 (s,x_1,x_2,x_3) = e^{ima_2(s)(x_1^2+x_2^2 )/2} \phi _2(x_1,x_2,x_3)$. By noting 
\begin{align*}
& (e^{i (t-s) \partial _3^2 /(2m)} e^{i\Omega (t-s) L_{(2)}} f)(x_1,x_2,x_3) \\ &= \frac{m}{(2\pi i (t-s))^{1/2}} \int_{\bf R} 
e^{im(x_3-z_3)^2/(t-s)} f(\hat{x}_1(t-s), \hat{x}_2 (t-s) , z_3) d z_3,
\end{align*}
where 
\begin{align*}
\MAT{\hat{x}_1 (t) \\ \hat{x}_2 (t)} = \MAT{\cos (\Omega (t)) & - \sin (\Omega (t)) \\ \sin (\Omega (t)) & \cos (\Omega (t))} \MAT{x _1 \\ x_2}, 
\end{align*}
(see, e.g., \S{3} of \cite{AK}), we have Lemma A.1. }

\end{document}